\documentclass[11pt]{article}
\usepackage{amsmath}
\usepackage{amsfonts}
\usepackage{amsthm}
\usepackage{verbatim}

\begin{document}

\def\N{\mathbb{N}}
\def\F{\mathbb{F}}
\def\Z{\mathbb{Z}}
\def\R{\mathbb{R}}
\def\Q{\mathbb{Q}}
\def\H{\mathcal{H}}

\parindent= 3.em \parskip=5pt

\centerline{\bf{On a particular hyperquadratic }} 
\centerline{\bf{ continued fraction in $\F(p)$ with $p>2$}}
\vskip 0.5 cm
\centerline{\bf{by A. Lasjaunias }}

\vskip 0.5 cm {\bf{Abstract.}} Given an odd prime number $p$, we describe a continued fraction in the field $\F(p)$ of power series in $1/T$ with coefficients in the finite field $\F_p$, where $T$ is a formal indeterminate. This continued fraction satisfies an algebraic equation of a particular type , with coefficients in $\F_p[T]$ which are explicitely given. We observe the close connection with other algebraic continued fractions which were studied thirty years ago by Mills and Robbins.
\vskip 0.5 cm
 In this note $p$ is an odd prime number and $\F_p$ is the finite field having $p$ elements. We denote $\F(p)$ the field of power series in $1/T$, where $T$ is a formal indeterminate, with coefficients in $\F_p$. Every irrational element $\alpha$ of $\F(p)$ is expanded as an infinite continued fraction denoted by $\alpha=[a_1,a_2,\cdots,a_n,\cdots$], where the $a_i$'s are polynomials in $\F_p[T]$. Our aim is to describe a particular sequence $(a_n)_{n\geq 1}$ such that $\alpha$ is algebraic over $\F_p(T)$ and it satisfies a particular equation explicitely described. This element $\alpha$ belongs to the subset of hyperquadratic elements : we have
$$A\alpha^{p+1}+B\alpha^p+C\alpha+D=0 \quad \text{ where } \quad A,B,C,D \in \F_p[T].$$
In the last forty years many examples of such hyperquadratic continued fractions were studied by different authors. For general information concerning this subject and also more references, the reader may consult \cite{L1} and \cite{S}. An important and fundamental article on continued fractions in function fields is due to Mills and Robbins \cite{MR}. The example described in the present note belongs to a family which has been derived from  this pioneer work. A detailed account about this connection was given in a recent note \cite{L2}.

 \par We recall briefly that if $\alpha=[a_1,a_2,\cdots,a_n,\cdots$], the tail of the expansion is denoted $\alpha_n=[a_n,\cdots]$, and we have 
$$\alpha=[a_1,a_2,\cdots,a_n,\alpha_{n+1}]=f_n(\alpha_{n+1})\quad \text{ for } \quad n\geq 1,$$
where $f_n$ is a linear fractional transformation with coefficients in $\F_p[T]$, depending on the first partial quotients $a_1,\cdots,a_n$. Indeed if we set $x_n/y_n=[a_1,\cdots,a_n]$ for the rational, called convergent, representing the finite continued fraction, we have
$$\alpha=(x_n\alpha_{n+1}+x_{n-1})/(y_n\alpha_{n+1}+y_{n-1})\quad \text{ for } \quad n\geq 1.\eqno{(1)}$$
Remember that $x_n$ and $y_n$ are the continuants obtained from the partial quotients by the same recurrence relation : $K_n=a_nK_{n-1}+K_{n-2}$ with the initial conditions $(x_1=a_1,x_0=1)$ for $x$ and $(y_1=1,y_0=0)$ for $y$.
\par Our continued fraction $\alpha$ is such that $\alpha=[u_1T,u_2T,u_3T,\alpha_4]$, where $(u_1,u_2,u_3)\in (\F_p^*)^3$ is an arbitrary triplet. Then we define a particular polynomial $F\in \F_p[T]$ by $F=(T^2+4)^{(p-1)/2}$. We will now build a  sequence in $\F_p[T]$, for each prime $p\geq 3$, depending on this triplet $(u_1,u_2,u_3)$ and this polynomial $F$.
\newline Let us define the following sequence $(P_n)_{n\geq 0}$ in $\F_p[T]$:
$$P_0=T \quad \text{ and } \quad P_{n+1}=F.P_n^p,\quad \text{ for } \quad n\geq 0.$$
Hence we have $\deg(P_{n+1})=p-1+p\deg(P_n)$. Consequently, we get $\deg(P_n)=2p^n-1$ for $n\geq 0$.
\newline For $n\geq 0$, we define two triples in $\F_p[T]$:
$$A_{1,n}=u_1T,u_2P_n,u_3T \quad \text{ and } \quad A_{2,n}=(4u_3)^{-1}T,4u_1u_2u_3P_n,(4u_1)^{-1}T .$$
For $n\geq 0$ we define two finite sequences of length $p^n-1$ in $\F_p[T]$:
$$B_{1,n}=((2u_3)^{-1}T,2u_3T)^{[(p^n-1)/2]} \quad \text{ and } \quad B_{2,n}=(2u_1T,(2u_1)^{-1}T)^{[(p^n-1)/2]},$$
where $a^{[k]}=a,\cdots,a\quad (k\quad \text{times})$. Finally, for $n\geq 0$, we define the following finite sequences in $\F_p[T]$:
$$C_{2n}=A_{1,2n}B_{1,2n}\quad \text{ and } \quad C_{2n+1}=A_{2,2n+1}B_{2,2n+1}.$$
 Note that $B_{1,0}=\emptyset$ and $C_0=u_1T,u_2T,u_3T$. We denote $R=T^p-TF \in \F_p[T]$. We observe that $R$ is simply the remainder in the Euclidean division of $T^p$ by $F$.
 Then we consider in $\F(p)$ the infinite continued fraction defined by
$$\alpha=[C_0,C_1,\cdots,C_n,\cdots]=[u_1T,u_2T,u_3T,(4u_3)^{-1}T,4u_1u_2u_3P_1,\cdots].$$ 
This continued fraction satisfies the following equalities $$\alpha=[u_1T,u_2T,u_3T,\alpha_4]\quad \text{ and } \quad \alpha^p=4u_1u_3F\alpha_4+u_1R.\eqno{(2)}$$
 \par Combining equalities $(1)$ and $(2)$, it becomes clear that $\alpha$ is algebraic and hyperquadratic.  Hence, we have
$$y_3\alpha^{p+1}-x_3\alpha^{p}+(4u_1u_3y_2F-u_1y_3R)\alpha+u_1x_3R-4u_1u_3x_2F=0.\eqno{(3)}$$

\par {\bf{\emph{Remark 1.}}} We shall now show the origin of the polynomial $F$ and its particular place in function fields arithmetic. In their article \cite[p.~400]{MR} Mills and Robbins introduce the following sequence $(f_n)_{n\geq 0}$ of polynomials in $\F_p[T]$:
$$f_0=1,\quad f_1=T\quad \text{ and } \quad f_{n}=Tf_{n-1}+f_{n-2}\quad \text{ for } \quad n\geq 2.$$
These polynomials can be regarded as the formal Fibonacci numbers. Indeed we have $\lim_n f_n/f_{n-1}=\omega=[T,T,\cdots,T,\cdots]$. This quadratic continued fraction $\omega$ is the analogue in the formal case of the real number $(1+\sqrt{5})/2=[1,1,\cdots,1,\cdots]$. These polynomials have been considered in several articles by the author (see \cite{L2} and the references therein). They satisfy, among others, the following identities
$$f_{p-1}=(T^2+4)^{(p-1)/2} \quad \text{ and } \quad f_p+f_{p-2}=T^p.$$
Note that $F=f_{p-1}$ and $R=T^p-TF=f_p+f_{p-2}-Tf_{p-1}=2f_{p-2}$. Mills and Robbins considered in $\F(p)$ for $p\geq 5$, the continued fractions defined by
$$\alpha=[u_1T,u_2T,\alpha_3]\quad \text{ and } \quad \alpha^p=F\alpha_3-(1/2)R,\eqno{(4)}$$
where $u_1\in \F_p$ with $u_1\neq 0,-1/2$ and $u_2=-u_1(1+2u_1)^{-1}\in \F_p^*$.
\newline In the same article they presented a general algorithm allowing to obtain the complete continued fraction expansion for certain hyperquadratic elements. They could apply this algorithm to obtain the explicit continued fraction of the elements satisfying (4). The remarkable fact about these continued fractions is that all the partial quotients are of the form $uT$ with $u\in \F_p^*$. If $u_1=-1$ then we have $\alpha=-\omega$, otherwise $\alpha$ is algebraic but not quadratic. This was the beginning of many extensions and generalizations concerning particular hyperquadratic continued fractions by the author and others (see \cite{L2} and the references therein).
\newline We have not tried to write a proof of our claim concerning the continued fraction defined by (2), using Mills and Robbins algorithm. However the truth of this claim is conforted by computer observations presented in the following remark.
\par {\bf{\emph{Remark 2.}}} Given an algebraic equation with coefficients in $\F_p[T]$ having a solution in fields of power series, there exists a mechanical process to obtain one after the other the partial quotients of this solution. This general principle was brought to light by Mohamed Mkaouar (Sfax university, Tunisia). With the help of Domingo Gomez (Cantabria university, Spain), using  Sage programming code, we could adapt this process to our equation (3). This program is presented below. By letting $p$ vary as well as the triplet $(u_1,u_2,u_3)$, the reader will be convinced of the correctness of the pattern for the continued fraction described above.
\begin{verbatim}

def contf(P,m):
    n = P.degree()
    an = P[n]
    an1 = P[n-1]
    a = []
    for i in range( m ):
        an = P[n]
        an1 = P[n-1]
        bar = - an1//an
        if P(bar) == 0:
            return bar
        else:
            P = P(x+bar)
            P=P.reverse()
            a.append(bar)
    return a

def leading_coefficients(P,m):
    return [f.leading_coefficient() for f in contf(P,m)]
def degree(P,m):
    return [f.degree() for f in contf(P,m)]


p = 7
F=GF(p)
u1=F(2)
u2=F(4)
u3=F(5)
z=F(4*u1*u3)
k=(p-1)/2
Ft.<t> = PolynomialRing(F)
Ftx.<x> = PolynomialRing(Ft)
Fib=(t^2+4)^k
R=t^p % Fib
x3=u1*u2*u3*t^3+(u1+u3)*t
y3=u2*u3*t^2+1
x2=u1*u2*t^2+1
y2=u2*t
P=y3*x^(p+1)-x3*x^p+(z*y2*Fib-u1*y3*R)*x+u1*x3*R-z*x2*Fib

print "p=",(p)
print "cfe",contf(P,7)
print "degrees",degree(P,65)
print "lead.coef.", leading_coefficients(P,65)


p= 7
cfe [2*t, 4*t, 5*t, 6*t, 6*t^13 + 2*t^11 + t^9 + 6*t^7, t, 4*t]
degrees [1, 1, 1, 1, 13, 1, 1, 1, 1, 1, 1, 1, 1, 97, 1, 1, 1, 1
, 1, 1, 1, 1, 1, 1, 1, 1, 1, 1, 1, 1, 1, 1, 1, 1, 1, 1, 1, 1, 1
, 1, 1, 1, 1, 1, 1, 1, 1, 1, 1, 1, 1, 1, 1, 1, 1, 1, 1, 1, 1, 1
, 1, 1, 1, 1, 685]
lead.coef. [2, 4, 5, 6, 6, 1, 4, 2, 4, 2, 4, 2, 2, 4, 5, 5, 3,
 5, 3, 5, 3, 5, 3, 5, 3, 5, 3, 5, 3, 5, 3, 5, 3, 5, 3, 5, 3, 5,
 3, 5, 3, 5, 3, 5, 3, 5, 3, 5, 3, 5, 3, 5, 3, 5, 3, 5, 3, 5, 3,
 5, 3, 5, 3, 6, 6]

\end{verbatim}

\par {\bf{\emph{Remark 3.}}} Contrarily to the historical example introduced by Mills and Robbins, in our continued fraction the sequence of the degrees of the partial quotients is unbounded. Indeed this sequence $(d_n)_{n\geq 1}$ contains the sub-sequence $(2p^n-1)_{n\geq 1}$. To describe more precisely this sequence of degrees, we introduce the sequence $(n_k)_{k\geq 1}$ such that $d_{n_k}=2p^k-1$ for $k\geq 1$. We have $n_1=5$. From the definition of the sequence of partial quotients we get
$$n_{k+1}=n_k+p^k+2 \quad \text{ and } \quad n_k=(p^k-1)/(p-1)+2k+2 \quad \text{ for } \quad k\geq 1.$$
Let us introduce $s_k=\sum_{1\leq i < n_k}d_i$ for $k\geq 1$. We have $s_1=4$. Since we have $d_n=1$ if $n\neq n_k$, we obtain, for $k\geq 1$,
$$s_{k+1}=s_k+2p^k-1+p^k+1=s_k+3p^k \quad \text{ and } \quad s_k=3(p^k-1)/(p-1)+1.$$
We are now interested in the irrationality  measure $\nu(\alpha)$ of our continued fraction $\alpha$. For any irrational $\alpha \in \F(p)$, if $\alpha=[a_1,a_2\cdots,a_n,\cdots]$, we have  
$$\nu(\alpha)=2+\limsup_n(\deg(a_{n+1})/\sum_{1\leq i\leq n}\deg(a_i)).$$ 
In our case, we can write
$$\limsup_n(d_{n+1}/\sum_{1\leq i\leq n}d_i)=\lim_k d_{n_k}/s_k.$$ 
Therefore, we obtain
 $$\nu(\alpha)=2+\lim_k(2p^k-1)/s_k=2+2(p-1)/3.$$ 
By a famous theorem of Liouville-Mahler, if $\alpha$ is algebraic over $\F_p(T)$ and $d=[\F_p(T,\alpha):\F_p(T)]>1$, then we have $\nu(\alpha)\in [2;d]$. In our case, since $\alpha$ satisfies the algebraic equation (3), this implies $2<\nu(\alpha)\leq d\leq p+1$. Hence in the particular case $p=3$, we must have $d=4$ and consequently, equation (3) is irreducible.

\vskip 0.5 cm
\begin{tabular}{ll}Alain LASJAUNIAS\\Institut de Math\'ematiques de Bordeaux  CNRS-UMR 5251
\\Universit\'e de Bordeaux \\Talence 33405, France \\E-mail: Alain.Lasjaunias@math.u-bordeaux.fr\\\end{tabular}

\end{document}